\documentclass[11pt]{article}
\usepackage[ansinew]{inputenc}
\usepackage{array}
\usepackage{color}
\usepackage{amsmath}
\usepackage{amsxtra}
\usepackage{amstext}
\usepackage{amssymb}
\usepackage{latexsym}
\usepackage{amsfonts,cite}
\usepackage{graphics,epstopdf}
\usepackage{epsfig}
\begin{document}

\sloppy
\newcommand{\proof}{{\it Proof~}}
\newtheorem{thm}{Theorem}[section]
\newtheorem{cor}[thm]{Corollary}
\newtheorem{lem}[thm]{Lemma}
\newtheorem{prop}[thm]{Proposition}
\newtheorem{eg}[thm]{Example}
\newtheorem{defn}[thm]{Definition}

\newtheorem{rem}[thm]{Remark}
\numberwithin{equation}{section}

\thispagestyle{empty}
\parindent=0mm

\begin{center}
{\Large{\bf Certain results for unified Apostol type-truncated exponential-Gould-Hopper polynomials and their relatives${}^{\star}$}}\\
\vspace{.25cm}

{\bf Serkan Araci${}^{\star,1}${\footnote{Corresponding author;E-mail:~serkan.araci@hku.edu.tr (Serkan Araci)}}, Mumtaz Riyasat${}^{2}$ {\footnote{E-mail:~mumtazrst@gmail.com (Mumtaz Riyasat)}},
Tabinda Nahid${}^{2}${\footnote{E-mail:~tabindanahid@gmail.com (Tabinda Nahid)}} and Subuhi Khan${}^{2}$ {\footnote{Email: subuhi2006@gmail.com (Subuhi Khan )}}}\\
\vspace{.25cm}

{\em $^{1}$Department of Economics, Faculty of Economics, Administrative and Social Sciences, Hasan Kalyoncu University, Gaziantep TR-27410, Turkey}\\
{\em $^{2}$Department of Mathematics, Faculty of Science, Aligarh Muslim University, Aligarh, India}\\
\end{center}

\begin{abstract}
\noindent
The present article aims to introduce a unified family of the Apostol type-truncated exponential-Gould-Hopper polynomials and to characterize
its properties via generating functions. A unified presentation of the generating function for the Apostol type-truncated exponential-Gould-Hopper
polynomials is established and its applications are given. By the use of operational techniques, the quasi-monomial properties for the unified family are proved. Several explicit representations and multiplication formulas related to these polynomials are obtained. Some general symmetric identities involving multiple power sums and Hurwitz-Lerch zeta functions are established by applying different analytical means on generating functions.

\end{abstract}
\vspace{.25cm}

\noindent
{\bf {\em Keywords:}} Unified Apostol type polynomials; Gould-Hopper polynomials; Truncated-exponential polynomials;
Explicit representations; Symmetric identities.\\

\noindent
{\em 2010 Mathematics Subject Classification:} 11B68; 05A10; 11B65.

\section{Introduction and preliminaries}

The generating functions have various applications in many branches of mathematics
and mathematical physics. These functions are defined by linear polynomials, differential
relations, globally referred to as functional equations. The functional equations arise in
well-defined combinatorial contexts and they lead systematically to well-defined classes
of functions. \\

The generating functions for the sequence of polynomials are used in analyzing
sequences of functions, in finding a closed formula for a sequence,
in finding recurrence relations and differential equations, relationships between sequences, asymptotic behavior of sequences,
in proving identities involving sequences and in solving enumeration problems in combinatorics and encoding their solutions.\\

On the subject of the Apostol-Bernoulli, Apostol-Euler and Apostol-Genocchi
polynomials and their various extensions, a remarkably large number of
investigations have appeared in the literature, see for example \cite{HE,Luo,Luo5,Luo2,Luo6,Luo9,Luo1}. Several authors obtained many interesting results involving various relatives of the Apostol-Bernoulli, Apostol-Euler and Apostol-Genocchi polynomials. The generalized Apostol Bernoulli polynomials $\mathfrak{B}_n^{(\alpha)}(x;\lambda)~(\alpha,~~\lambda \in \mathbb{C};~~x \in \mathbb{R})$ are defined by the generating function of the form \cite{Luo1}:
\begin{equation}
\left(\frac{t}{\lambda e^{t}-1}\right)^\alpha e^{xt}=\sum\limits_{n=0}^{\infty} \mathfrak{B}_n^{(\alpha)}(x;\lambda)\frac{t^{n}}{n!},~~~~~~~|t| < 2\pi, ~\textrm{when}~\lambda=1;~|t|<|\textrm{log} \lambda|,~\textrm{when}~\lambda\neq1.
\end{equation}

The Apostol-Bernoulli polynomials and
Apostol-Bernoulli numbers can be obtained from the generalized Apostol-Bernoulli polynomials by
\begin{equation}
\mathfrak{B}_n^{(1)}(x;1)=\mathfrak{B}_n(x;\lambda)
\hspace{1cm}\mathfrak{B}_n(0;\lambda)=\mathfrak{B}_n(\lambda),~~~~~n \in \mathbb{N}_0,
\end{equation}
respectively. The case $\lambda=1$ in the above relations gives the classical Bernoulli polynomials $B_n(x)$ and Bernoulli numbers
$B_n$ \cite{Erd3}.\\

The generalized Apostol-Euler polynomials $\mathfrak{E}_n^{(\alpha)}(x;\lambda)~(\alpha,~~\lambda \in \mathbb{C};~~x \in \mathbb{R})$ are defined by the generating function of the form \cite{Luo}:
\begin{equation}
\left(\frac{2}{\lambda e^{t}+1}\right)^\alpha e^{xt}=\sum\limits_{n=0}^{\infty} \mathfrak{E}_n^{(\alpha)}(x;\lambda)\frac{t^{n}}{n!},~~~~~~~|t| < \pi, ~\textrm{when}~\lambda=1;~|t|<|\textrm{log} (-\lambda)|,~\textrm{when}~\lambda\neq1.
\end{equation}

The Apostol-Euler polynomials and
Apostol-Euler numbers can be obtained from the generalized Apostol-Euler polynomials by
\begin{equation}
\mathfrak{E}_n^{(1)}(x;1)=\mathfrak{E}_n(x;\lambda)
\hspace{1cm}\mathfrak{E}_n(0;\lambda)=\mathfrak{E}_n(\lambda),~~~~~n \in \mathbb{N}_0,
\end{equation}
respectively. The case $\lambda=1$ in the above relations gives the classical Euler polynomials $E_n(x)$ and Euler numbers
$E_n$ \cite{Erd3}.\\

The generalized Apostol-Genocchi polynomials $\mathcal{G}_n^{(\alpha)}(x;\lambda)~(\alpha,~~\lambda \in \mathbb{C};~~x \in \mathbb{R})$ are defined by the generating function of the form \cite{Luo2}:
\begin{equation}
\left(\frac{2t}{\lambda e^{t}+1}\right)^\alpha e^{xt}=\sum\limits_{n=0}^{\infty} \mathcal{G}_n^{(\alpha)}(x;\lambda)\frac{t^{n}}{n!},~~~~~~|t| < \pi, ~\textrm{when}~\lambda=1;~|t|<|\textrm{log} (-\lambda)|,~\textrm{when}~\lambda\neq1.
\end{equation}

The Apostol-Genocchi polynomials and
Apostol-Genocchi numbers can be obtained from the generalized Apostol-Genocchi polynomials by
\begin{equation}
\mathcal{G}_n^{(1)}(x;1)=\mathcal{G}_n(x;\lambda)
\hspace{1cm}\mathcal{G}_n(0;\lambda)=\mathcal{G}_n(\lambda),~~~~~n \in \mathbb{N}_0,
\end{equation}
respectively. The case $\lambda=1$ in the above relations gives the classical Genocchi polynomials $G_n(x)$ and Genocchi numbers
$G_n$ \cite{Sandor}.\\

Various unified forms of the Apostol-Bernoulli, Apostol-Euler and Apostol-Genocchi polynomials are introduced and studied by many authors in a systematic manner
via different analytic means and generating functions method, see for example \cite{Ozarlanunified,Kurtsymmetry,Kurtsymmetry1}. We recall the following unified form of the Apostol-Bernoulli, Apostol-Euler and Apostol-Genocchi polynomials \cite{Ozarlanunified}:\\

For $\alpha,~\beta \in \mathbb{C}$, $a,~b\in \mathbb{R}\setminus \{0\}$ and $k\in \mathbb{N}_{0},$ the
\emph{unified Apostol type polynomials} $\mathcal{P}_{n,\beta
}^{(\alpha )}(x;k,a,b)$  are defined by the following
generating function:
\begin{equation}\label{1.7}
\left( \frac{2^{1-k}t^{k}}{\beta ^{b}~e^t-a^{b}}\right) ^{\alpha
}e^{xt}=\sum\limits_{n=0}^{\infty }\mathcal{P}_{n,\beta}^{(\alpha )}(x;k,a,b)\frac{t^{n}}{n!}.
\end{equation}

In fact, the following special cases hold:
\begin{equation}\label{1.8}
\mathcal{P}_{n,\lambda}^{(\alpha)}(x;1,1,1):=\mathfrak{B}
_{n}^{(\alpha)}(x;\lambda),
\end{equation}
\begin{equation}
\mathcal{P}_{n,\lambda}^{(\alpha)}(x;0,-1,1):=\mathfrak{E}
_{n}^{(\alpha)}(x;\lambda),
\end{equation}
\begin{equation}
\mathcal{P}_{n,\frac{\lambda}{2}}^{(\alpha)}\left(x;1,-\frac{1}{2},1\right):=
\mathcal{G}_{n}^{(\alpha)}(x;\lambda),
\end{equation}
where $\mathfrak{B}_{n}^{(\alpha)}(x;\lambda)$, $\mathfrak{E}
_{n}^{(\alpha)}(x;\lambda)$ and $\mathcal{G}_{n}^{(\alpha)}(x;
\lambda)$ are the generalized forms of the Apostol-Bernoulli, Apostol-Euler and
Apostol-Genocchi polynomials.\\

The multi-variable forms of the special functions serve as an analytical foundation for the majority of problems in mathematical physics that have been solved exactly and find broad practical applications. These polynomials have provided new means of analysis for the solution of large classes of partial differential equations often encountered in physical problems. We recall certain 2-variable special polynomials of mathematical physics.\\

The Gould-Hopper polynomials (GHP) $H_n^{(m)}(x,y)$ are defined by means of the following generating function \cite[p.58 (6.3)]{Gould}:
\begin{equation}\label{1.11}
e^{xt+yt^m}=\sum\limits_{n=0}^\infty H_n^{(m)}(x,y)\frac{t^n}{n!},~~~~~m>2
\end{equation}
and possess the following series expansion:
\begin{equation}
H_n^{(m)}(x,y)=n!\sum\limits_{k=0}^{[\frac{n}{m}]}\frac{y^k~x^{n-mk}}{k!~(n-mk)!}.
\end{equation}

The 2-variable Hermite Kamp$\acute{e}$ de Feriet polynomials (2VHKdFP) $H_n(x,y)$ \cite{App1} can be obtained from the GHP $H_n^{(m)}(x,y)$
such that $H_n(x,y):=H_n^{(2)}(x,y)$,
which on choosing $x=2x,~~y=-1$ gives the classical Hermite polynomials $H_n(x)$ \cite{Andrew}.\\

The 2-variable truncated exponential polynomials (2VTEP) (of order $r$) $e_n^{(r)}(x,z)$ are defined by means of the following
generating function \cite[p.174 (30)]{Bessel}:
\begin{equation}\label{1.13}
\frac{e^{xt}}{(1-zt^{r})}=\sum\limits_{n=0}^{\infty }e_{n}^{(r)}(x,z)\frac{t^{n}}{n!}
\end{equation}
and possess the following series expansion:
\begin{equation}
e_{n}^{(r)}(x,z)=n!\sum_{k=0}^{[\frac{n}{r}]}~\frac{z^k~x^{n-rk}}{(n-rk)!}.
\end{equation}

The 2-variable truncated exponential polynomials (2VTEP) \cite{Trunc} can be obtained from the 2VTEP (of order $r$) $e_n^{(r)}(x,z)$, such that
$e_{n}^{(2)}(x,z)=n!{}_{[2]}e_{n}(x,z)$, which on choosing $z=1$ gives the classical truncated exponential polynomials $n!{}_{_{[2]}}e_n(x)$ \cite{Andrew}.\\

To introduce the hybridized forms of the unified Apostol type polynomials and to characterize their properties via generating functions will lead to a new approach.
The process of combining two multi variable forms of special polynomials either by replacement technique or operational technique
is called hybridization of polynomials, which shows that the properties of new hybridized polynomial lie within the context of parent polynomials. Several investigations have made to introduce and study hybrid families of special polynomials systematically via various analytic means
\cite{2VATP,Khansymmetry}.\\

Here, we first introduce the 3-variable truncated exponential-Gould-Hopper polynomials (3VTEGHP) ${}_eH_{n}^{(m,r)}(x,y,z)$ by means of the following generating function:
\begin{equation}\label{1.15}
\frac{e^{xt+yt^m}}{(1-zt^{r})}=\sum\limits_{n=0}^{\infty }{}_eH_{n}^{(m,r)}(x,y,z)\frac{t^n}{n!}.
\end{equation}

{\em {\bf Lines of proof}: Simply by expanding the exponential function and on replacing the powers of $x$, that is $x^n~(n=0,1,2,\ldots)$ by the polynomials $e_n^{(r)}(x,z)~(n=0,1,2,\ldots)$ in generating function (1.11) of the GHP $H_n^{(m)}(x,y)$ and on use of generating equation (1.13) of the 2VTEP (of order $r$) $e_n^{(r)}(x,z)$ we get desired equation \eqref{1.15}.}\\

We find the following series expansion for the 3VTEGHP ${}_eH_n^{(m,r)}(x,y,z)$:
\begin{equation}
{}_eH_n^{(m,r)}(x,y,z)=n!\sum\limits_{k=0}^{[\frac{n}{r}]}\sum\limits_{l=0}^{[\frac{n-rk}{m}]}\frac{z^k~x^{n-rk-ml}y^l}{l!~(n-rk-ml)!}.
\end{equation}

Using similar approach and by taking 3VTEGHP ${}_eH_{n}^{(m,r)}(x,y,z)$ given by \eqref{1.15} as base in generating equation (1.7) of the unified Apostol type polynomials, we find the following generating function for the unified Apostol type-truncated exponential-Gould-Hopper polynomials (UATyTEGHP ) ${_{{}_eH}}\mathcal{P}_{n,\beta}^{(\alpha,m,r)}(x,y,z;k,a,b)$:
\begin{equation}\label{1.16}
\left( \frac{2^{1-k}t^{k}}{\beta ^{b}~e^t-a^{b}}\right) ^{\alpha
}\frac{e^{xt+yt^m}}{(1-zt^{r})}=\sum\limits_{n=0}^{\infty}{_{{}_eH}}\mathcal{P}_{n,\beta}^{(\alpha,m,r)}(x,y,z;k,a,b)\frac{t^n}{n!}.
\end{equation}

Using equations (1.7), (1.15) and (1.16) in equation (1.17), we get the following series expansion for the UATyTEGHP ${_{{}_eH}}\mathcal{P}_{n,\beta}^{(\alpha,m,r)}(x,y,z;k,a,b)$:
\begin{equation}
{_{{}_eH}}\mathcal{P}_{n,\beta}^{(\alpha,m,r)}(x,y,z;k,a,b)=n!\sum\limits_{p=0}^n \sum\limits_{k=0}^{[\frac{p}{r}]}\sum\limits_{l=0}^{[\frac{p-rk}{m}]}\frac{\mathcal{P}_{n-p}^{(\alpha)}(k,a,b)z^k~x^{p-rk-ml}y^l}{(n-p)!l!~(p-rk-ml)!}.
\end{equation}

For $y=0$, the UATyTEGHP reduce to the truncated-exponential based Apostol-type polynomials considered in \cite{aaaa} and for $z=0$ becomes the Gould-Hopper based Apostol-type polynomials \cite{bbbb}. It should be noted that the UATyTEGHP ${_{{}_eH}}\mathcal{P}_{n,\beta}^{(\alpha,m,r)}(x,y,z;k,a,b)$ include the following special cases:\\

\noindent
{\bf Remark 1.1.} Setting $r=m=2$, equation \eqref{1.16} yields the following generating function for the 3-variable unified Apostol type-truncated exponential-Hermite polynomials (3VUATyTEHP ) ${_{{}_{[2]}eH}}\mathcal{P}_{n,\beta}^{(\alpha)}(x,y,z;k,a,b)$:
\begin{equation}
\left( \frac{2^{1-k}t^{k}}{\beta ^{b}~e^t-a^{b}}\right) ^{\alpha
}\frac{e^{xt+yt^2}}{(1-zt^{2})}=\sum\limits_{n=0}^{\infty}{_{{}_{[2]}eH}}\mathcal{P}_{n,\beta}^{(\alpha)}(x,y,z;k,a,b)t^n.
\end{equation}

For $x \rightarrow 2x,~y=-1;~z=1$, equation (1.19) yields the generating function for the unified Apostol type-truncated exponential-Hermite polynomials (UATyTEHP ) ${_{{}_{[2]}eH}}\mathcal{P}_{n,\beta}^{(\alpha)}(x;k,a,b)$:
\begin{equation}
\left( \frac{2^{1-k}t^{k}}{\beta ^{b}~e^t-a^{b}}\right) ^{\alpha
}\frac{e^{2xt-t^2}}{(1-t^{2})}=\sum\limits_{n=0}^{\infty}{_{{}_{[2]}eH}}\mathcal{P}_{n,\beta}^{(\alpha)}(x;k,a,b)t^n.
\end{equation}

\noindent
{\bf Remark 1.2.} Setting $k = a = b = 1$ and $\beta=\lambda$ and using ${_{{}_eH}}\mathcal{P}_{n,\lambda}^{(\alpha,m,r)}(x,y,z;1,1,1):={_{{}_eH}}\mathfrak{B}_{n}^{(\alpha,m,r)}(x,y,z;\lambda)$ in equation \eqref{1.16}, we get the truncated exponential-Gould-Hopper-Apostol Bernoulli polynomials (TEGHABP) ${_{{}_eH}}\mathfrak{B}_{n}^{(\alpha,m,r)}(x,y,z;\lambda)$ given by
\begin{equation}\label{1.21}
\left( \frac{t}{\lambda e^t-1}\right) ^{\alpha
}\frac{e^{xt+yt^m}}{(1-zt^{r})}=\sum\limits_{n=0}^{\infty }{_{{}_eH}}\mathfrak{B}_{n}^{(\alpha,m,r)}(x,y,z;\lambda)\frac{t^n}{n!}.
\end{equation}

\noindent
{\bf Remark 1.3.} Setting $k =0,~ a =-1,~ b = 1$ and $\beta=\lambda$ and using ${_{{}_eH}}\mathcal{P}_{n,\lambda}^{(\alpha,m,r)}(x,y,z;0,-1,1):={_{{}_eH}}\mathfrak{E}_{n}^{(\alpha,m,r)}(x,y,z;\lambda)$ in equation \eqref{1.16}, we get the truncated exponential-Gould-Hopper-Apostol Euler polynomials (TEGHAEP)${_{{}_eH}}\mathfrak{E}_{n}^{(\alpha,m,r)}(x,y,z;\lambda)$ given by
\begin{equation}\label{1.22}
\left( \frac{2}{\lambda e^t+1}\right) ^{\alpha
}\frac{e^{xt+yt^m}}{(1-zt^{r})}=\sum\limits_{n=0}^{\infty }{_{{}_eH}}\mathfrak{E}_{n}^{(\alpha,m,r)}(x,y,z;\lambda)\frac{t^n}{n!}.
\end{equation}

\noindent
{\bf Remark 1.4.} Setting $k =1,~ a =-1/2,~ b =1$ and $\beta=\lambda/2$ and using ${_{{}_eH}}\mathcal{P}_{n,\lambda/2}^{(\alpha,m,r)}(x,y,z;1,-1/2,1):={_{{}_eH}}\mathcal{G}_{n}^{(\alpha,m,r)}(x,y,z;\lambda)$ in equation \eqref{1.16}, we get the truncated exponential-Gould-Hopper-Apostol Genocchi polynomials (TEGHAGP) ${_{{}_eH}}\mathcal{G}_{n}^{(\alpha,m,r)}(x,y,z;\lambda)$ given by
\begin{equation}\label{1.23}
\left( \frac{2t}{\lambda e^t+1}\right) ^{\alpha
}\frac{e^{xt+yt^m}}{(1-zt^{r})}=\sum\limits_{n=0}^{\infty }{_{{}_eH}}\mathcal{G}_{n}^{(\alpha,m,r)}(x,y,z;\lambda)\frac{t^n}{n!}.
\end{equation}

This paper aims at presenting the study of multi-variable special polynomials which play an important role in several fields of physics, applied mathematics and engineering. The present article is an attempt to further stress the importance of the truncated exponential-Gould-Hopper and unified Apostol type polynomials to introduce a hybrid family of the unified Apostol type-truncated exponential-Gould-Hopper polynomials. In Section 2, the quasi-monomial properties of these polynomials are derived via operational techniques. In Sections 3 and 4, several explicit representations and multiplication formulas are obtained. In Section 5, some general symmetry identities are established.

\section{Quasi-monomial properties}

We know that the operational techniques including differential operators provide a systematic and analytic approach to study special functions.
Differential equations have been the primary motivation for the introduction of these techniques. The operational techniques combined with the monomiality principle \cite{Dat} open new possibilities to deal with the theoretical foundations of special polynomials and also to introduce new families of special polynomials.\\

According to the monomiality principle, there exist two operators $\Phi^+$ and $\Phi^-$ playing, respectively, the role of multiplicative and derivative operators for a polynomial set $\{p_n(x)\}_{n \in \mathbb{N}}$, that is, $\Phi^+$ and $\Phi^-$ satisfy the following identities, for all $n \in \mathbb{N}$:
\begin{equation}\label{2.1}
\Phi^{+}\{p_n(x)\}=p_{n+1}(x),\hspace{1cm}
\Phi^{-}\{p_n(x)\}=n~p_{n-1}(x).
\end{equation}

The polynomial set $\{p_n(x)\}_{n \in \mathbb{N}}$ is then called a quasi-monomial. If $\Phi^+$ and $\Phi^-$ have differential realizations, then the polynomial set $\{p_n(x)\}_{n \in \mathbb{N}}$ satisfy the following differential equation:
\begin{equation}\label{2.2}
\Phi^{+}\Phi^{-}\{p_n(x)\}=n~p_n(x).
\end{equation}

In order to establish the quasi-monomiality of the
UATyTEGHP ${_{{}_eH}}\mathcal{P}_{n,\beta}^{(\alpha,m,r)}(x,y,z;k,a,b)$, we prove the following results:

\begin{thm}
The unified Apostol type-truncated exponential-Gould-Hopper polynomials ${_{{}_eH}}\mathcal{P}_{n,\beta}^{(\alpha,m,r)}(x,y,z;k,a,b)$ are quasi-monomial with respect to the following multiplicative and derivative operators:
\begin{equation}\label{aa}
\begin{array}{lllll}
& \hspace{1.5cm}\Phi^+_{{}_{eH}\mathcal{P}^{(\alpha,m,r)}}:=  x + my~ ( \partial_x)^{m-1} - \frac{rz ~( \partial_x)^{r-1}}{(1-z( \partial_x)^r)} + \frac{\alpha k}{ \partial_x} - \frac{\alpha ~\beta^b e^{ \partial_x}}{\beta^b e^{ \partial_x} - a^b},
\end{array}
\end{equation}
\begin{equation}\label{b}
\Phi^-_{{}_{eH}\mathcal{P}^{(\alpha,m,r)}}:= \partial_x. \hspace{5.29cm}
\end{equation}
\end{thm}

\noindent
\begin{proof}
We consider the following identity:
\begin{equation}\label{c}
\partial_x \left\{ \left( \frac{2^{1-k} t^k}{\beta^b e^t - a^b}\right)^\alpha  \frac{e^{xt + yt^m}}{(1- z t^r)}\right\} = t \left( \frac{2^{1-k} t^k}{\beta^b e^t - a^b}\right)^\alpha  \frac{e^{xt + yt^m}}{(1- z t^r)}.
\end{equation}

Differentiating both sides of the equation \eqref{1.16} w.r.t. $t$, it follows that
\begin{equation}\label{d}
\begin{split}
\Bigg( x + my~  t^{m-1} - \frac{rz t^{r-1}}{(1-z t^r)} &+ \frac{\alpha k}{t} - \frac{\alpha ~\beta^b e^t}{\beta^b e^t - a^b} \Bigg) \Bigg(\frac{2^{1-k} t^k}{\beta^b e^t - a^b} \Bigg) \frac{e^{xt +yt^m}}{(1- zt^r)} \\
&= \sum\limits_{n=0}^\infty n~ {_{{}_eH}}\mathcal{P}_{n+1,\beta}^{(\alpha,m,r)}(x, y, z; k,a,b)\frac{t^{n-1}}{n !},
\end{split}
\end{equation}
which in view of identity \eqref{c} becomes
\begin{equation}\label{e}
\begin{split}
\Big( x + my~ ( \partial_x)^{m-1} - \frac{rz ~( \partial_x)^{r-1}}{(1-z( \partial_x)^r)} + \frac{\alpha k}{ \partial_x} - \frac{\alpha ~\beta^b e^{ \partial_x}}{\beta^b e^{ \partial_x} - a^b} \Big) & \sum\limits_{n=0}^\infty  {_{{}_eH}}\mathcal{P}_{n,\beta}^{(\alpha,m,r)}(x, y, z; k,a,b) \frac{t^n}{n!}\\
& \hspace{-3cm}= \sum\limits_{n=0}^\infty  {_{{}_eH}}\mathcal{P}_{n+1,\beta}^{(\alpha,m,r)}(x, y, z; k,a,b) \frac{t^n}{n!}.
\end{split}
\end{equation}

On equating the coefficients of same power of $t$ in both sides of the above equation and in view of monomiality principle equation (2.1), assertion \eqref{aa} is proved.\\

Use of equation \eqref{1.16} in identity \eqref{c} yields
\begin{equation}\label{f}
 \partial_x  \left\{ \sum\limits_{n=0}^\infty  {_{{}_eH}}\mathcal{P}_{n,\beta}^{(\alpha,m,r)}(x, y, z; k,a,b)~ \frac{t^n}{n!} \right\}= \sum\limits_{n=0}^\infty n~ {_{{}_eH}}\mathcal{P}_{n-1,\beta}^{(\alpha,m,r)}(x, y, z; k,a,b)~ \frac{t^n}{n!}.
\end{equation}

On equating the coefficients of same power of $t$ in both sides of the above equation and in view of monomiality principle equation (2.1), assertion \eqref{b} follows.\\
\end{proof}

\noindent
{\bf Remark 2.1.} By making use of  multiplicative and derivative operators given by expressions \eqref{aa} and \eqref{b} in equation (2.2), we deduce the following consequence of Theorem 2.1:

\begin{cor}
The unified Apostol type-truncated exponential-Gould-Hopper polynomials ${_{{}_eH}}\mathcal{P}_{n,\beta}^{(\alpha,m,r)}(x,y,z;k,a,b)$ satisfy the following differential equation:
\begin{equation}\label{g}
\begin{split}
\Big( x + my ( \partial_x)^{m} + \alpha k - \frac{rz ( \partial_x)^{r}}{(1-z( \partial_x)^r)} - \frac{\alpha ~\beta^b e^{\partial_x}  \partial_x}{\beta^b e^{ \partial_x} - a^b}  - n\Big)  {_{{}_eH}}\mathcal{P}_{n,\beta}^{(\alpha,m,r)}(x, y, z; k,a,b) =0.
\end{split}
\end{equation}
\end{cor}

In view of Remarks 1.2-1.4, we find the differential operators and differential equations for the TEGHABP ${_{{}_eH}}\mathfrak{B}_{n,\beta}^{(\alpha,m,r)}(x,y,z;\lambda)$, TEGHAEP ${_{{}_eH}}\mathfrak{E}_{n,\beta}^{(\alpha,m,r)}(x,y,z;\lambda)$
and TEGHAGP ${_{{}_eH}}\mathcal{G}_{n,\beta}^{(\alpha,m,r)}(x,y,z;\lambda)$. We present these results in Table 1:\\

\noindent
\textbf{Table 1.~~Results for the TEGHABP ${_{{}_eH}}\mathfrak{B}_{n,\beta}^{(\alpha,m,r)}(x,y,z;k,a,b)$, TEGHAEP}

\textbf{\hspace{.7cm}${_{{}_eH}}\mathfrak{E}_{n,\beta}^{(\alpha,m,r)}(x,y,z;\lambda)$ and  TEGHAGP ${_{{}_eH}}\mathcal{G}_{n,\beta}^{(\alpha,m,r)}(x,y,z;\lambda)$}\\
\\
{\tiny{
\begin{tabular}{lllll}
\hline
&&&&\\
{\bf S.}  & {\bf Special}  &{\bf Differential} & {\bf Differential}\\
{\bf No.}&{\bf polynomials}& {\bf operators}& {\bf equations}\\
\hline
&&&\\
{\bf I.}  & {\bf TEGHABP} & $\Phi^+_{{}_{eH}\mathfrak{B}^{(\alpha,m,r)}}:= \Big( x + my~ ( \partial_x)^{m-1}$& $\Big( x + my (\partial_x)^{m} + \alpha - \frac{rz ( \partial_x)^{r}}{(1-z(\partial_x)^r)} - \frac{\alpha ~\lambda e^{ \partial_x}  \partial_x}{\lambda e^{ \partial_x} - 1} - n \Big)$\\
&${_{{}_eH}}\mathfrak{B}_{n,\beta}^{(\alpha,m,r)}(x,y,z;k,a,b)$ & $- \frac{rz ~( \partial_x)^{r-1}}{(1-z( \partial_x)^r)} + \frac{\alpha }{ \partial_x} - \frac{\alpha ~\lambda e^{ \partial_x}}{\lambda e^{ \partial_x} - 1} \Big)$  &$   {_{{}_eH}}\mathfrak{B}_{n}^{(\alpha,m,r)}(x, y, z;\lambda) =0$\\
&&$\Phi^-_{{}_{eH}\mathfrak{B}^{(\alpha,m,r)}}:=\partial_x$&\\
\hline
&&&\\
{\bf II.}  & {\bf TEGHAEP} & $\Phi^+_{{}_{eH}\mathfrak{E}^{(\alpha,m,r)}}:= \Big( x + my~ ( \partial_x)^{m-1}$& $\Big( x + my ( \partial_x)^{m} - \frac{rz ( \partial_x)^{r}}{(1-z( \partial_x)^r)} - \frac{\alpha ~\lambda e^{ \partial_x}  \partial_x}{\lambda e^{ \partial_x} +1} -n \Big)$\\
&${_{{}_eH}}\mathfrak{E}_{n,\beta}^{(\alpha,m,r)}(x,y,z;\lambda)$& $- \frac{rz ~( \partial_x)^{r-1}}{(1-z( \partial_x)^r)} - \frac{\alpha ~\lambda e^{ \partial_x}}{\lambda e^{ \partial_x}+ 1} \Big)$  &$  {_{{}_eH}}\mathfrak{E}_{n}^{(\alpha,m,r)}(x, y, z;\lambda) =0$\\
&&$\Phi^-_{{}_{eH}\mathfrak{E}^{(\alpha,m,r)}}:= \partial_x$&\\
\hline
&&&\\
{\bf III.}  & {\bf TEGHAGP} & $\Phi^+_{{}_{eH}\mathcal{G}^{(\alpha,m,r)}}:=  \Big( x + my~ ( \partial_x)^{m-1}$& $\Big( x + my ( \partial_x)^{m} + \alpha - \frac{rz ( \partial_x)^{r}}{(1-z( \partial_x)^r)} - \frac{\alpha ~\lambda e^{ \partial_x}  \partial_x}{2(\lambda e^{ \partial_x}+ 1)} - n \Big)$\\
&${_{{}_eH}}\mathcal{G}_{n,\beta}^{(\alpha,m,r)}(x,y,z;\lambda)$& $- \frac{rz ~( \partial_x)^{r-1}}{(1-z( \partial_x)^r)} + \frac{\alpha }{ \partial_x} - \frac{\alpha ~\lambda/2 e^{ \partial_x}}{2(\lambda e^{ \partial_x} +1)} \Big)$  &$ {_{{}_eH}}\mathcal{G}_{n}^{(\alpha,m,r)}(x, y, z;\lambda) =0$\\
&&$\Phi^-_{{}_{eH}\mathcal{G}^{(\alpha,m,r)}}:= \partial_x$&\\
\hline
\end{tabular}}}\\
\vspace{.35cm}

In the next section, we establish certain explicit representations for the UATyTEGHP ${_{{}_eH}}\mathcal{P}_{n,\beta}^{(\alpha,m,r)}(x,y,z;k,a,b)$ and for their special cases.

\section{Explicit representations}

In order to derive the explicit representations for the UATyTEGHP ${_{{}_eH}}\mathcal{P}_{n,\beta}^{(\alpha,m,r)}(x,y,z;k,a,b)$, we recall the following definitions:\\

\noindent
{\bf Definition 3.1.} Let $x \in \mathbb{R}$, $\lambda\in \mathbb{C}$ and $\nu \in \mathbb{N}_0$. The array type polynomials $S(x;n,\nu,\lambda)$ are defined by \cite{Simsek1}:
\begin{equation}
\frac{(\lambda e^{t}-1)^{\nu}}{\nu!}~e^{xt}=\sum\limits_{n=0}^\infty S(x;n,\nu;\lambda)\frac{t^n}{n!}.
\end{equation}
where ${S}(n,v;\lambda):=S(0;n,\nu,\lambda)$ are $\lambda$-Stirling type numbers
of second kind, which again for $\lambda=1$ becomes the Stirling numbers of second kind $S(n,k)$ \cite{Comet}.\\

\noindent
{\bf Definition 3.2.} The generalized Hurwitz-Lerch Zeta function $\Phi_\mu(z,s,a)$ \cite{Goyal} is defined by
\begin{equation}
 \Phi_\mu(z,s,a)=\sum\limits_{n=0}^\infty \frac{(\mu)_n}{n!}~\frac{z^n}{(n+a)^s},
\end{equation}
which for $\mu=1$ becomes the Hurwitz-Lerch Zeta function $\Phi(z,s,a)$ \cite{Hwtz}.\\

Now, we derive explicit representations for the UATyTEGHP ${_{{}_eH}}\mathcal{P}_{n,\beta}^{(\alpha,m,r)}(x,y,z;k,a,b)$ by proving the following results:

\begin{thm}
 The following explicit formula for the unified Apostol type-truncated exponential-Gould-Hopper polynomials $ {_{{}_eH}}\mathcal{P}_{n,\beta}^{(\alpha,m,r)}(x,y,z;k,a,b)$ in terms of the array type polynomials $S(x;n,\nu;\lambda)$ holds true:
\begin{equation}
 {_{{}_eH}}\mathcal{P}_{n-k\alpha,\beta}^{(-\alpha,m,r)}(x,y,z;k,a,b)
=\frac{(n-k\alpha)!\alpha!a^{b\alpha}}{2^{(1-k)\alpha}}~\sum\limits_{p=0}^{[\frac{n}{m}]}\sum\limits_{s=0}^{[\frac{n}{r}]}
\frac{S(x;n-rs-mp,\nu;(\frac{\beta}{a})^b )y^p z^s}{(n-rs-mp)!p!},
\end{equation}
$\hspace{8.5cm}a,~b \in \mathbb{R}^+$, $x,~y,~z \in \mathbb{R}$; $\lambda \in \mathbb{C}$;~~~ $n \in \mathbb{N}_0.$
\end{thm}

\noindent
\begin{proof}
Rewriting generating equation \eqref{1.16} in the following form:
\begin{equation}
\sum\limits_{n=0}^\infty {_{{}_eH}}\mathcal{P}_{n,\beta}^{(-\alpha,m,r)}(x,y,z;k,a,b)\frac{t^n}{n!}
=\left( \frac{2^{1-k}t^{k}}{\beta ^{b}~e^t-a^b}\right) ^{-\alpha}\frac{e^{xt+yt^m}}{(1-zt^{r})},
\end{equation}
which on simplification becomes
\begin{equation}
\sum\limits_{n=0}^\infty {_{{}_eH}}\mathcal{P}_{n,\beta}^{(-\alpha,m,r)}(x,y,z;k,a,b)\frac{t^{n+k\alpha}}{n!}=\frac{\alpha!a^{b\alpha}}{2^{(1-k)\alpha}}~\left(\frac{((\frac{\beta}{a})^b ~e^t-1)^\alpha}{\alpha!}~e^{xt}\right)e^{yt^m}(1-zt^{r})^{-1}.
\end{equation}

Using equation (3.1) and making appropriate expansions in equation (3.5) and then using the Cauchy product rule in the resultant equation, we are led to assertion (3.3).\\
\end{proof}

Similarly, we can prove the following result:
\begin{thm}
The following explicit formula for the unified Apostol type-truncated exponential-Gould-Hopper polynomials ${_{{}_eH}}\mathcal{P}_{n,\beta}^{(\alpha,m,r)}(x,y,z;k,a,b)$ in terms of the 3-variable truncated exponential-Gould-Hopper polynomials ${_e}H_n^{(m,r)}(x,y,z)$ holds true:
\begin{equation}
 {_{{}_eH}}\mathcal{P}_{n-k\alpha,\beta}^{(-\alpha,m,r)}(x,y,z;k,a,b)
=\frac{\alpha!a^{b\alpha}}{2^{(1-k)\alpha}}\frac{(n-k\alpha)!}{n !}~\sum\limits_{p=0}^n {n \choose p}
S\Big(p,\alpha;\Big(\frac{\beta}{a}\Big)^b\Big){}_eH_{n-p}^{(m,r)}(x,y,z).
\end{equation}
\end{thm}

\vspace{.1cm}
\begin{thm}
The following explicit formula for the unified Apostol type-truncated exponential-Gould-Hopper polynomials ${_{{}_eH}}\mathcal{P}_{n,\beta}^{(\alpha,m,r)}(x,y,z;k,a,b)$ in terms of the generalized Hurwitz-Lerch Zeta function  $\Phi_\mu(z,s,a)$ holds true:

\begin{equation}
{_{{}_eH}}\mathcal{P}_{n+k\alpha,\beta}^{(\alpha,m,r)}(x,y,z;k,a,b)= \frac{(n+k\alpha)!}{n!} \left(\frac{2^{1-k}}{a^b}\right)^\alpha\sum\limits_{l=0}^n {n \choose l}
\Phi_\alpha\left(\Big(\frac{\beta}{a}\Big)\Big)^{b},-l,x\right) {}_eH_{n-l}^{(m,r)}(0,y,z),
\end{equation}
\hspace{8cm} $~~a,~b \in \mathbb{R}^+$, $x,~y,~z \in \mathbb{R}$; $\lambda \in \mathbb{C}$;~~ $n \in \mathbb{N}_0$.
\end{thm}

\noindent
\begin{proof}
Considering generating function (1.17) and simplifying it follows that
\begin{equation}
\sum\limits_{n=0}^\infty {_{{}_eH}}\mathcal{P}_{n,\beta}^{(\alpha,m,r)}(x,y,z;k,a,b)\frac{t^{n-k\alpha}}{n!}=\left(\frac{2^{1-k}}{a^b}\right)^\alpha\sum\limits_{k=0}^\infty
(\alpha)_k~\Big(\frac{\beta}{a}\Big)^{bk} \frac{e^{(k+x)t}}{k!}e^{yt^m}(1-zt^{r})^{-1}.
\end{equation}

In view of equations (1.15) and (3.2), we have
\begin{equation}
\sum\limits_{n=0}^\infty {_{{}_eH}}\mathcal{P}_{n,\beta}^{(\alpha,m,r)}(x,y,z;k,a,b)\frac{t^{n-k\alpha}}{n!}=\left(\frac{2^{1-k}}{a^b}\right)^\alpha\sum\limits_{l=0}^\infty
\Phi_\alpha\left(\Big(\frac{\beta}{a}\Big)\Big)^{b},-l,x\right) \frac{t^l}{l!} \sum\limits_{n=0}^\infty {}_eH_{n}^{(m,r)}(0,y,z) \frac{t^n}{n!}.
\end{equation}

Using Cauchy product rule in the above equation and then equating the coefficients of same powers of $t$ in both sides of resultant equation, assertion (3.7) is proved.\\
\end{proof}

Next, we derive an implicit formula for the UATyTEGHP ${_{{}_eH}}\mathcal{P}_{n,\beta}^{(\alpha,m,r)}(x,y,z;k,a,b)$ by proving the following result:
\begin{thm}
The following implicit formula for the unified Apostol type-truncated exponential-Gould-Hopper polynomials  ${_{{}_eH}}\mathcal{P}_{n,\beta}^{(\alpha,m,r)}(x,y,z;k,a,b)$ holds true:
\begin{equation}
 {_{{}_eH}}\mathcal{P}_{n-k\gamma,\beta}^{(\alpha-\gamma,m,r)}(x,y,z;k,a,b)= \frac{(n-k\gamma)! \gamma!}{n!} \left(\frac{a^b}{2^{1-k}}\right)^\gamma
 \sum\limits_{n=0}^l {n \choose l}{_{{}_eH}}\mathcal{P}_{n-l,\beta}^{(\alpha,m,r)}(x,y,z;k,a,b) S\Big(l,\gamma,\Big(\frac{\beta}{a}\Big)^b\Big).
\end{equation}
\end{thm}

\noindent
\begin{proof}
Rewriting generating equation \eqref{1.16} in the following form:
\begin{equation}
\sum\limits_{n=0}^\infty {_{{}_eH}}\mathcal{P}_{n-k \gamma,\beta}^{(\alpha - \gamma,m,r)}(x,y,z;k,a,b) \frac{t^{n-k\gamma}}{(n-k \gamma)!}
=\left( \frac{2^{1-k}t^{k}}{\beta ^{b}~e^t-a^b}\right) ^{\alpha - \gamma}\frac{e^{xt+yt^m}}{(1-zt^{r})},
\end{equation}
which on simplification becomes
\begin{equation}
\begin{split}
\sum\limits_{n=0}^\infty {_{{}_eH}}\mathcal{P}_{n -k \gamma,\beta}^{(\alpha -\gamma,m,r)}(x,y,z;k,a,b)\frac{t^{n}}{(n- k \gamma)!}& =\frac{\gamma!a^{b \gamma}}{2^{(1-k)\gamma}}~\left(\frac{((\frac{\beta}{a})^b ~e^t-1)^\gamma}{\gamma !}\right)\\
& \times \sum\limits_{n=0}^\infty {_{{}_eH}}\mathcal{P}_{n,\beta}^{(\alpha,m,r)}(x,y,z;k,a,b)\frac{t^{n}}{n!}.
\end{split}
\end{equation}

Using equation (3.1) and making appropriate expansions in equation (3.12) and then using the Cauchy product rule in the resultant equation, we are led to assertion (3.10).\\
\end{proof}

In view of Remarks 1.2-1.4, we find the explicit representations for the TEGHABP ${_{{}_eH}}\mathfrak{B}_{n,\beta}^{(\alpha,m,r)}(x,y,z;\lambda)$, TEGHAEP ${_{{}_eH}}\mathfrak{E}_{n,\beta}^{(\alpha,m,r)}(x,y,z;\lambda)$
and TEGHAGP ${_{{}_eH}}\mathcal{G}_{n,\beta}^{(\alpha,m,r)}(x,y,z;\lambda)$. We present these results in Table 2:\\

\noindent
\textbf{Table 2.~~Explicit representations for the  TEGHABP ${_{{}_eH}}\mathfrak{B}_{n,\beta}^{(\alpha,m,r)}(x,y,z;\lambda)$, TEGHAEP}

\textbf{\hspace{.7cm} ${_{{}_eH}}\mathfrak{E}_{n,\beta}^{(\alpha,m,r)}(x,y,z;\lambda)$
and TEGHAGP ${_{{}_eH}}\mathcal{G}_{n,\beta}^{(\alpha,m,r)}(x,y,z;\lambda)$} \\
\\
{\tiny{
\begin{tabular}{lll}
\hline
&&\\
{\bf S.No.}  & {\bf Special polynomials} &{\bf Explicit representations}\\
\hline
&&\\
{\bf I.}  & {\bf TEGHABP} & $ {_{{}_eH}}\mathfrak{B}_{n-\alpha}^{(-\alpha,m,r)}(x,y,z;\lambda)
={(n-\alpha)!\alpha!}~\sum\limits_{p=0}^{[n/m]}\sum\limits_{s=0}^{[n/r]}
\frac{S(x;n-rs-mp,\nu;\lambda)y^p z^s}{(n-rs-mp)!p!}$\\
&${_{{}_eH}}\mathfrak{B}_{n,\beta}^{(\alpha,m,r)}(x,y,z;\lambda)$& $ {_{{}_eH}}\mathfrak{B}_{n-\alpha}^{(-\alpha,m,r)}(x,y,z;\lambda)
=\frac{\alpha!}{n!}(n-\alpha)!~\sum\limits_{p=0}^n {n \choose p}
S\big(p,\alpha;\lambda\big){}_eH_{n-p}^{(m,r)}(x,y,z)$\\
&& ${_{{}_eH}}\mathfrak{B}_{n-\gamma}^{(\alpha-\gamma,m,r)}(x,y,z;\lambda)= \frac{(n-\gamma)! \gamma!}{n!}
 \sum\limits_{n=0}^l {n \choose l}{_{{}_eH}}\mathfrak{B}_{n-l}^{(\alpha,m,r)}(x,y,z;\lambda) S\big(l,\gamma,\lambda\big)$\\
 && ${_{{}_eH}}\mathfrak{B}_{n+\alpha}^{(\alpha,m,r)}(x,y,z;\lambda)= \frac{(n+\alpha)!}{n!}\sum\limits_{l=0}^n {n \choose l}
\Phi_\alpha\left(\lambda,-l,x\right) {}_eH_{n-l}^{(m,r)}(0,y,z)$\\
&&\\
\hline
&&\\
{\bf II.}  & {\bf TEGHAEP} & $ {_{{}_eH}}\mathfrak{E}_{n}^{(-\alpha,m,r)}(x,y,z;\lambda)
=\frac{n!\alpha!(-1)^\alpha}{2^{k\alpha}}~\sum\limits_{p=0}^{[n/m]}\sum\limits_{s=0}^{[n/r]}
\frac{S(x;n-rs-mp,\nu;-\lambda)y^p z^s}{(n-rs-mp)!p!}$\\
&${_{{}_eH}}\mathfrak{E}_{n,\beta}^{(\alpha,m,r)}(x,y,z;\lambda)$& $ {_{{}_eH}}\mathfrak{E}_{n}^{(-\alpha,m,r)}(x,y,z;\lambda)
=\frac{\alpha!(-1)^\alpha}{2^\alpha}~\sum\limits_{p=0}^n {n \choose p}
S\big(p,\alpha;-\lambda\big){}_eH_{n-p}^{(m,r)}(x,y,z)$\\
&& ${_{{}_eH}}\mathfrak{E}_{n}^{(\alpha-\gamma,m,r)}(x,y,z;\lambda)=\gamma!(2)^\gamma
 \sum\limits_{n=0}^l {n \choose l}{_{{}_eH}}\mathfrak{E}_{n-l}^{(\alpha,m,r)}(x,y,z;\lambda) S\big(l,\gamma,-\lambda\big)$\\
 && ${_{{}_eH}}\mathfrak{E}_{n}^{(\alpha,m,r)}(x,y,z;\lambda)=(-2)^\alpha\sum\limits_{l=0}^n {n \choose l}
\Phi_\alpha\left(-\lambda,-l,x\right) {}_eH_{n-l}^{(m,r)}(0,y,z)$\\
&&\\
\hline
&&\\
{\bf III.}  & {\bf TEGHAGP} & ${_{{}_eH}}\mathcal{G}_{n-\alpha}^{(-\alpha,m,r)}(x,y,z;\lambda)
={(n-\alpha)!\alpha!(-1/2)^\alpha}~\sum\limits_{p=0}^{[n/m]}\sum\limits_{s=0}^{[n/r]}
\frac{S(x;n-rs-mp,\nu;-\lambda)y^p z^s}{(n-rs-mp)!p!}$\\
&${_{{}_eH}}\mathcal{G}_{n,\beta}^{(\alpha,m,r)}(x,y,z;\lambda)$& ${_{{}_eH}}\mathcal{G}_{n-\alpha}^{(-\alpha,m,r)}(x,y,z;\lambda)
=\frac{\alpha!}{n!}(-1/2)^\alpha(n-\alpha)!~\sum\limits_{p=0}^n {n \choose p}
S\big(p,\alpha;-\lambda\big){}_eH_{n-p}^{(m,r)}(x,y,z)$\\
&& ${_{{}_eH}}\mathcal{G}_{n}^{(\alpha-\gamma,m,r)}(x,y,z;\lambda)=\frac{(n-\gamma)!}{n!}\gamma!(-1/2)^\gamma
 \sum\limits_{n=0}^l {n \choose l}{_{{}_eH}}\mathcal{G}_{n-l}^{(\alpha,m,r)}(x,y,z;\lambda) S\big(l,\gamma,-\lambda\big)$\\
 && ${_{{}_eH}}\mathcal{G}_{n+\alpha}^{(\alpha,m,r)}(x,y,z;\lambda)=\frac{(n+\alpha)!}{n!}(-2)^\alpha\sum\limits_{l=0}^n {n \choose l}
\Phi_\alpha\left(-\lambda,-l,x\right) {}_eH_{n-l}^{(m,r)}(0,y,z)$\\
&&\\
\hline
\end{tabular}}}\\
\vspace{.5cm}

In the next section, we derive multiplication formulae for the UATyTEGHP ${_{{}_eH}}\mathcal{P}_{n,\beta}^{(\alpha,m,r)}(x,y,z;k,a,b)$ and for their special cases.\\

\section{Multiplication formulae}

The multiplication formulas related to the Apostol-Bernoulli, Apostol-Euler and Apostol-Genocchi polynomials and for their unifications are derived in \cite{Luo5,Luo6,Luo9}. In this section, we derive the multiplication formula for the UATyTEGHP ${_{{}_eH}}\mathcal{P}_{n,\beta}^{(\alpha,m,r)}(x,y,z;k,a,b)$. The multiplication formulas for the TEGHABP ${_{{}_eH}}\mathfrak{B}_{n,\beta}^{(\alpha,m,r)}(x,y,z;\lambda)$, TEGHAEP ${_{{}_eH}}\mathfrak{E}_{n,\beta}^{(\alpha,m,r)}(x,y,z;\lambda)$
and TEGHAGP ${_{{}_eH}}\mathcal{G}_{n,\beta}^{(\alpha,m,r)}(x,y,z;\lambda)$ are deduced as special cases. We recall the following definitions:\\

\noindent
{\bf Definition 3.1.} ({\em Multinomial identity \cite{Comet}}) If $x_1,x_2,\ldots,x_m$ are commuting
elements of a ring ($\Leftrightarrow x_i x_j=x_jx_i,~~~1\leq i<j \leq m$), then we have for all integers $\alpha \geq0$:
\begin{equation}\label{2.7}
(1+x_1+x_2+\ldots+x_m)^\alpha=\sum\limits_{\nu_1,\nu_2,\ldots,\nu_m \geq 0}{\alpha \choose \nu_1,\nu_2,\ldots,\nu_m}x_1^{\nu_1}x_2^{\nu_2}\ldots x_m^{\nu_m},
\end{equation}
the last summation takes over all positive or zero integers $a_i \geq 0$ such that $a_1 + a_2 + a_m = n$,
\begin{equation}
{\alpha \choose \nu_1,\nu_2,\ldots,\nu_{m}}=\frac{\alpha!}{\nu_1! \nu_2! \ldots \nu_m!},
\end{equation}
are called multinomial coefficients.\\

\noindent
{\bf Definition 3.2.} ({\em Generalized Multinomial identity \cite{Comet}}) If $x_1,x_2,\ldots,x_m$ are commuting
elements of a ring ($\Leftrightarrow x_i x_j=x_jx_i,~~~1\leq i<j \leq m$), then we have for all integers real or complex $\alpha$:
\begin{equation}
(x_1+x_2+\ldots+x_m)^\alpha=\sum\limits_{\substack{0\leq \nu_1\leq \nu_2\leq \nu_{m}\leq \alpha\\\nu_1+\nu_2+\ldots+\nu_{m}=\alpha}} {\alpha \choose \nu_1,\nu_2,\ldots,\nu_m}x_1^{\nu_1}x_2^{\nu_2}\ldots x_m^{\nu_m},
\end{equation}
the last summation takes place over all positive or zero integers $\nu_i \geq 0$,
\begin{equation}
{\alpha \choose \nu_1,\nu_2,\ldots,\nu_{m}}=\frac{[\alpha]_{\nu_1+\nu_2+\ldots+\nu_{m}}}{\nu_1! \nu_2! \ldots \nu_m!}=\frac{\alpha(\alpha-1)(\alpha-2)\ldots{(\alpha-\nu_1-\nu_2-\nu_m+1)}}{\nu_1! \nu_2! \ldots \nu_m!},
\end{equation}
are called generalized multinomial coefficients.\\

In order to find multiplication formula for the UATyTEGHP ${_{{}_eH}}\mathcal{P}_{n,\beta}^{(\alpha,m,r)}(x,y,z;k,a,b)$, the following result is proved.

\begin{thm}
Let $a,~b \in \mathbb{R}^+$, $x,~y,~z \in \mathbb{R}$; $k \in \mathbb{C}$; $s \in \mathbb{N}$ and $n \in \mathbb{N}_0$. The following multiplication formula for the unified Apostol type-truncated exponential-Gould-Hopper polynomials ${_{{}_eH}}\mathcal{P}_{n,\beta}^{(\alpha,m,r)}(x,y,z;k,a,b)$ holds true:
\begin{equation}\label{2.4}
\begin{array}{lllll}
{_{{}_eH}}\mathcal{P}_{n,\beta}^{(\alpha,m,r)}(sx,s^my,s^rz;k,a,b)& = s^{n-k\alpha}~ a^{(s-1)b\alpha} \sum\limits_{\nu_1,\nu_2,\ldots,\nu_{s-1} \geq 0}{\alpha \choose \nu_1,\nu_2,\ldots,\nu_{s-1}}\\
& \times \left(\frac{\beta}{a} \right)^{bl}~{_{{}_eH}}\mathcal{P}_{n,\beta^s}^{(\alpha,m,r)}(x+\frac{l}{s},y,z;k,a^s,b^s),~~~~for ~~~s ~~~~odd,
\end{array}
\end{equation}

\begin{equation}\label{2.5}
\begin{array}{lllll}
{_{{}_eH}}\mathcal{P}_{n,\beta}^{(p,m,r)}(sx,s^my,s^rz;k,a,b)& = \frac{(-1)^p~ 2^{(1-k) p} s^{n-kp}}{(n+1)_{(1-k)p}}~ a^{(s-1)bp} \sum\limits_{\substack{0\leq \nu_1\leq \nu_2\leq \nu_{s-1}\leq p\\\nu_1+\nu_2+\ldots+\nu_{s-1}=p}}{p \choose \nu_1,\nu_2,\ldots,\nu_{s-1}}\\
& \times \left(\frac{\beta}{a} \right)^{bl}~{_{{}_eH}}\mathfrak{B}_{{n+(1-k)p},\beta^s}^{(p,m,r)}(x+\frac{l}{s},y,z;k,a^s,b^s),~~~~for ~~~s ~~~~even,
\end{array}
\end{equation}
where $\nu_1+2\nu_2+\ldots+(s-1)\nu_{s-1} =l$.
\end{thm}

\noindent
\noindent
\begin{proof}
We consider generating equation (1.17) in the following form:
\begin{equation}\label{2.5}
\sum\limits_{n=0}^\infty {_{{}_eH}}\mathcal{P}_{n,\beta}^{(\alpha,m,r)}(sx,s^my,s^rz;k,a,b) \frac{t^n}{n!} = \frac{1}{s^{k\alpha}}\left(\frac{2^{1-k} (st)^k}{\beta^{bs} e^{ts} - a^{bs}}\right)^\alpha \left( \frac{\beta^{bs} e^{ts} - a^{bs}}{\beta^b e^t -a^b}\right)^\alpha \frac{e^{xst+ y(st)^m}}{(1-s^rz(t)^r)},
\end{equation}
which on simplification gives
\begin{equation}\label{2.6}
\begin{array}{lllll}
\sum\limits_{n=0}^\infty {_{{}_eH}}\mathcal{P}_{n,\beta}^{(\alpha,m,r)}(sx,s^my,s^rz;k,a,b)\frac{t^n}{n!}& = \frac{1}{s^{k\alpha}}\left(\frac{2^{1-k} (st)^k}{\beta^{bs} e^{ts} - a^{bs}}\right)^\alpha a^{(s-1)b\alpha} \left( \frac{ -\left(\frac{\beta}{a}\right)^{bs} e^{ts} +1}{-\left(\frac{\beta}{a}\right)^{b} e^{t} +1}\right)^\alpha \frac{e^{xst+ y(st)^m}}{(1-z(st)^r)}\\
& \hspace{-1cm}= \frac{1}{s^{k\alpha}}\left(\frac{2^{1-k} (st)^k}{\beta^{bs} e^{ts} - a^{bs}}\right)^\alpha a^{(s-1)b\alpha} \left(\sum\limits_{h=0}^{s-1} \left( \left(\frac{\beta}{a} \right)^b e^t \right)^h \right)^\alpha \frac{e^{xst+ y(st)^m}}{(1-z(st)^r)}.
\end{array}
\end{equation}

On use of equation (4.1), the above equation becomes
\begin{equation}\label{2.8}
\begin{array}{lllll}
\sum\limits_{n=0}^\infty {_{{}_eH}}\mathcal{P}_{n,\beta}^{(\alpha,m,r)}(sx,s^my,s^rz;k,a,b) \frac{t^n}{n!} &= \frac{1}{s^{k\alpha}}\left(\frac{2^{1-k} (st)^k}{\beta^{bs} e^{ts} - a^{bs}}\right)^\alpha a^{(s-1)b\alpha} \sum\limits_{\nu_1,\nu_2,\ldots,\nu_{s-1} \geq 0}{\alpha \choose \nu_1,\nu_2,\ldots,\nu_{s-1}}\\
& \hspace{-1cm} \times \left(  \left(\frac{\beta}{a} \right)^b\right)^{\nu_1+2\nu_2+\ldots+(s-1)\nu_{s-1}} e^{(\nu_1+2\nu_2+\ldots+(s-1)\nu_{s-1})t}~ \frac{e^{xst+ y(st)^m}}{(1-z(st)^r)},
\end{array}
\end{equation}
or, equivalently
\begin{equation}\label{2.8}
\begin{array}{lllll}
\sum\limits_{n=0}^\infty {_{{}_eH}}\mathcal{P}_{n,\beta}^{(\alpha,m,r)}(sx,s^my,s^rz;k,a,b) \frac{t^n}{n!}& = \frac{1}{s^{k\alpha}}\left(\frac{2^{1-k} (st)^k}{\beta^{bs} e^{ts} - a^{bs}}\right)^\alpha a^{(s-1)b\alpha} \sum\limits_{\nu_1,\nu_2,\ldots,\nu_{s-1} \geq 0}{\alpha \choose \nu_1,\nu_2,\ldots,\nu_{s-1}}\\
& \times  \left(\frac{\beta}{a} \right)^{bl} ~ \frac{e^{\left(x+ \frac{l}{s}\right)st+ y(st)^m}}{(1-z(st)^r)}.
\end{array}
\end{equation}

Again, using equation (1.17) in the above equation, we have
\begin{equation}\label{2.9}
\begin{array}{lllll}
\sum\limits_{n=0}^\infty {_{{}_eH}}\mathcal{P}_{n,\beta}^{(\alpha,m,r)}(sx,s^my,s^rz;k,a,b)\frac{t^n}{n!} &= \sum\limits_{n=0}^\infty \Big(s^{n-k\alpha} a^{(s-1)b\alpha} \sum\limits_{\nu_1,\nu_2,\ldots,\nu_{s-1} \geq 0}{\alpha \choose \nu_1,\nu_2,\ldots,\nu_{s-1}}\\
& \times  \left(\frac{\beta}{a} \right)^{bl} {_{{}_eH}}\mathcal{P}_{n,\beta^s}^{(\alpha,m,r)}\left(x+\frac{l}{s},y,z;k,a^s,b^s\right)\Big)\frac{t^n}{n!}.
\end{array}
\end{equation}

Equating the coefficients of identical power of $t$ on both sides of the above equation,  assertion \eqref{2.4} is proved. With use of equation (4.2)
(for $r$ even), we can similarly prove assertion (4.6).\\
\end{proof}

In view of Remarks 1.2-1.4, we find the multiplication formulas for the TEGHABP ${_{{}_eH}}\mathfrak{B}_{n,\beta}^{(\alpha,m,r)}(x,y,z;\lambda)$, TEGHAEP ${_{{}_eH}}\mathfrak{E}_{n,\beta}^{(\alpha,m,r)}(x,y,z;\lambda)$
and TEGHAGP ${_{{}_eH}}\mathcal{G}_{n,\beta}^{(\alpha,m,r)}(x,y,z;\lambda)$. We present these results in Table 3:\\

\noindent
\textbf{Table 3.~~Multiplication formulas for the TEGHABP ${_{{}_eH}}\mathfrak{B}_{n,\beta}^{(\alpha,m,r)}(x,y,z;\lambda)$, TEGHAEP }

\textbf{ \hspace{.7cm}${_{{}_eH}}\mathfrak{E}_{n,\beta}^{(\alpha,m,r)}(x,y,z;\lambda)$
and TEGHAGP ${_{{}_eH}}\mathcal{G}_{n,\beta}^{(\alpha,m,r)}(x,y,z;\lambda)$ } \\
\\
{\tiny{
\begin{tabular}{lll}
\hline
&&\\
{\bf S.No.}  & {\bf Special polynomials} &{\bf Multiplication formulas}\\
\hline
&&\\
{\bf I.}  & {\bf TEGHABP} & ${_{{}_eH}}\mathfrak{B}_{n}^{(\alpha,m,r)}(sx,s^my,s^r z;\lambda) = s^{n - \alpha} (\lambda)^l$\\
& ${_{{}_eH}}\mathfrak{B}_{n,\beta}^{(\alpha,m,r)}(x,y,z;\lambda)$& $\sum\limits_{\nu_1,\nu_2,\ldots,\nu_{s-1} \geq 0} {\alpha \choose \nu_1,\nu_2,\ldots,\nu_{s-1}} {_{{}_eH}}\mathfrak{B}_{n}^{(\alpha,m,r)}(x+\frac{l}{s},y,z;\lambda^s) $\\
&&\\
\hline
&&\\
{\bf II.}  & {\bf TEGHAEP} & ${_{{}_eH}}\mathfrak{E}_{n}^{(\alpha,m,r)}(sx,s^my,s^r z;\lambda) = s^{n}  (-1)^{(s-1)\alpha}(-\lambda)^l$\\
& ${_{{}_eH}}\mathfrak{E}_{n,\beta}^{(\alpha,m,r)}(x,y,z;\lambda)$& $\sum\limits_{\nu_1,\nu_2,\ldots,\nu_{s-1} \geq 0} {\alpha \choose \nu_1,\nu_2,\ldots,\nu_{s-1}} {_{{}_eH}}\mathfrak{E}_{n}^{(\alpha,m,r)}(x+\frac{l}{s},y,z;\lambda^s) $\\
&&${_{{}_eH}}\mathfrak{E}_{n}^{(p,m,r)}(sx,s^my,s^rz;\lambda)= \frac{(-2)^{p} s^{n}}{(n+1)_{p}}~ (-1)^{(s-1)p} \sum\limits_{\substack{0\leq \nu_1\leq \nu_2\leq \nu_{s-1}\leq p\\\nu_1+\nu_2+\ldots+\nu_{s-1}=p}}{p \choose \nu_1,\nu_2,\ldots,\nu_{s-1}}$\\
& &$\times (-\lambda)^{l}~{_{{}_eH}}\mathfrak{B}_{{n+p}}^{(p,m,r)}(x+\frac{l}{s},y,z;\lambda^s)$\\
&&\\
\hline
&&\\
{\bf III.}  & {\bf TEGHAGP}& ${_{{}_eH}}\mathcal{G}_{n}^{(\alpha,m,r)}(sx,s^my,s^r z;\lambda) = s^{n - \alpha} \left(-\frac{1}{2}\right)^{(s-1)\alpha} \left(-\lambda\right)^l$\\
&${_{{}_eH}}\mathcal{G}_{n,\beta}^{(\alpha,m,r)}(x,y,z;\lambda)$ &$\sum\limits_{\nu_1,\nu_2,\ldots,\nu_{s-1} \geq 0} {\alpha \choose \nu_1,\nu_2,\ldots,\nu_{s-1}} {_{{}_eH}}\mathcal{G}_{n}^{(\alpha,m,r)}(x+\frac{l}{s},y,z;\lambda^s) $\\
&&${_{{}_eH}}\mathcal{G}_{n,\beta}^{(p,m,r)}(sx,s^my,s^rz;k,a,b) = (-1)^p~ s^{n-p} \Big(\frac{-1}{2}\Big)^{(s-1)p} \sum\limits_{\substack{0\leq \nu_1\leq \nu_2\leq \nu_{s-1}\leq p\\\nu_1+\nu_2+\ldots+\nu_{s-1}=p}}{p \choose \nu_1,\nu_2,\ldots,\nu_{s-1}}$\\
& &$\times (-\lambda)^{l}~{_{{}_eH}}\mathfrak{B}_{{n}}^{(p,m,r)}(x+\frac{l}{s},y,z;\lambda^s)$\\
&&\\
\hline
\end{tabular}}}\\
\vspace{.35cm}

In the next section, we derive certain symmetry identities for the UATyTEGHP ${_{{}_eH}}\mathcal{P}_{n,\beta}^{(\alpha,m,r)}(x,y,z;k,a,b)$  and for their special cases.

\section{Symmetric identities}

The symmetry identities related to several extensions of the Apostol type polynomials are established by many authors \cite{Kurtsymmetry,Kurtsymmetry1,Khansymmetry}. This provides motivation to establish symmetry identities for the unified Apostol type-truncated exponential-Gould-Hopper polynomials ${_{{}_eH}}\mathcal{P}_{n,\beta}^{(\alpha,m,r)}(x, y,z; k,a,b)$ . First, we review the following definitions:\\

\noindent
{\bf Definition~5.1.}~For an arbitrary real or complex parameter $\lambda$, a number $\mathcal{S}_k(n;\lambda)$ is defined by the
following generating function:
\begin{equation}\label{3.1}
\frac{\lambda e^{(n+1)t}-1}{\lambda e^{t}-1}=\sum\limits_{k=0}^{\infty} \mathcal{S}_{k}(n;\lambda)\frac{t^k}{k!}.
\end{equation}

\noindent
{\bf Definition~5.2.}~For an arbitrary real or complex parameter $\lambda$, the multiple power sums $\mathcal{S}_{k}^{\left( l\right) }\left( m;\lambda \right)$ is defined by the
following generating function:
\begin{equation}\label{3.2}
\left( \frac{1-\lambda ^{m}e^{mt}}{1-\lambda e^{t}}\right) ^{l}=\frac{1}{\lambda ^{l}}%
\sum\limits_{n=0}^{\infty }\left\{ \sum\limits_{p=0}^{n}\binom{n}{p}\left(
-l\right) ^{n-p}\mathcal{S}_{k}^{\left( l\right) }\left( m;\lambda \right) \right\}
\frac{t^{n}}{n!}.
\end{equation}

In order to derive the symmetry identities for the UATyTEGHP ${_{{}_eH}}\mathcal{P}_{n,\beta}^{(\alpha,m,r)}(x,y,z;k,a,b)$, the following results are proved:

\begin{thm} For all integers $l,~q >0$ and $n \geq0, ~\alpha \geq1, ~a,b,k \in \mathbb{C}$, the following symmetry relation for the unified Apostol type-truncated exponential-Gould-Hopper polynomials ${_{{}_eH}}\mathcal{P}_{n,\beta}^{(\alpha,m,r)}(x,y,z;k,a,b)$  holds true:
\begin{equation}\label{3.4}
\begin{array}{lll}
& \hspace{.5cm}\sum\limits_{j=0}^{n} l^{n-j} q^{j+1} {_{{}_eH}}\mathcal{P}_{n-j,\beta}^{(\alpha,m,r)}(qx,q^my,q^rz;k,a,b) \sum\limits_{i=0}^j S_i\left(l-1, \left( \frac{\beta}{a} \right)^b\right)\\
& \hspace{1.5cm} \times{_{{}_eH}}\mathcal{P}_{j-i,\beta}^{(\alpha-1,m,r)}(lX,l^mY,l^rZ;k,a,b)\\
\\
& \hspace{.5cm}=\sum\limits_{j=0}^{n} q^{n-j} l^{j+1} {_{{}_eH}}\mathcal{P}_{n-j,\beta}^{(\alpha,m,r)}(lx,l^my,l^rz;k,a,b) \sum\limits_{i=0}^j S_i\left(q-1, \left( \frac{\beta}{a} \right)^b\right)\\
&\hspace{1.5cm} \times {_{{}_eH}}\mathcal{P}_{j-i,\beta}^{(\alpha-1,m,r)}(qX,q^mY,q^rZ;k,a,b).
\end{array}
\end{equation}
\end{thm}

\noindent
\begin{proof}Let
\begin{equation}\label{3.5}
G(t):= \frac{2^{(1-k)(2\alpha -1)} t^{k(2\alpha -1)} e^{lqxt + y(lqt)^m} (\beta^b e^{lqt} -a^b) e^{lqXt + Y(lqt)^m}}{(\beta^b e^{lt}- a^b)^\alpha~ (\beta^b e^{qt}- a^b)^\alpha~ (1-z(lqt)^r)~ (1 - Z(lqt)^r)},
\end{equation}
which on rearranging the powers becomes
\begin{equation}\label{3.6}
G(t) = \frac{1}{l^{k\alpha} q^{k(\alpha -1)}} \left(\frac{2^{1-k} (lt)^k}{\beta^b e^{lt} - a^b} \right)^\alpha \frac{e^{qx(lt) + q^m y(lt)^m}}{(1-q^r z(lt)^r)} \left( \frac{\beta^b e^{lqt} -a^b}{\beta^b e^{qt} -a^b} \right) \left(\frac{2^{1-k} (qt)^k }{\beta^b e^{qt} - a^b} \right)^{\alpha -1} \frac{e^{lX(qt) + l^mY(qt)^m}}{(1 - l^rZ(qt)^r)}.
\end{equation}

Using equations (1.17) and \eqref{3.1} in r.h.s. of equation \eqref{3.6} and applying the Cauchy product rule, we find
\begin{equation}\label{3.8}
\begin{array}{lll}
G(t) =  &\frac{1}{l^{k\alpha} q^{k(\alpha -1)}} \sum\limits_{n=0}^\infty \Bigg( \sum\limits_{j=0}^{n}\frac{l^{n-j} q^j}{i !}{_{{}_eH}}\mathcal{P}_{n-j,\beta}^{(\alpha,m,r)}(qx, q^my, q^rz;k,a,b) \sum\limits_{i=0}^{j}
\mathcal{S}_{i}\left(l-1;\left(\frac{\beta}{a}\right)^b \right)\\
& \times {_{{}_eH}}\mathcal{P}_{j-i,\beta}^{(\alpha -1,m,r)}(lX,l^mY, l^rZ;k,a,b) \Bigg) t^n.
\end{array}
\end{equation}
Similarly, we have
\begin{equation}\label{3.9}
\begin{array}{lll}
G(t) =  &\frac{1}{q^{k\alpha} l^{k(\alpha -1)}} \sum\limits_{n=0}^\infty \Bigg( \sum\limits_{j=0}^{n}\frac{q^{n-j} l^j}{i !}{_{{}_eH}}\mathcal{P}_{n-j,\beta}^{(\alpha,m,r)}(lx, l^my, l^rz;k,a,b) \sum\limits_{i=0}^{j}
\mathcal{S}_{i}\left(q-1;\left(\frac{\beta}{a}\right)^b \right)\\
& \times {_{{}_eH}}\mathcal{P}_{j-i,\beta}^{(\alpha -1,m,r)}(qX,q^mY, q^rZ;k,a,b) \Bigg) t^n.
\end{array}
\end{equation}

Equating the coefficients of same powers of $t$ in r.h.s. of expansions \eqref{3.8} and \eqref{3.9}, we are led to assertion \eqref{3.4}.
\end{proof}

\begin{thm} For each pair of positive integers $c,~d$ and for all integers $n \geq0, ~\alpha \geq1, ~a,b,k \in \mathbb{C}$, the following symmetry relation for the unified Apostol type-truncated exponential-Gould-Hopper polynomials ${_{{}_eH}}\mathcal{P}_{n,\beta}^{(\alpha,m,r)}(x,y,z;k,a,b)$  holds true:
\begin{equation}\label{3.10}
\begin{array}{llll}
&\hspace{-.86cm}\sum\limits_{l=0}^n c^{n-l} d^l \sum\limits_{i=0}^{c-1}\sum\limits_{j=0}^{d-1}~\left(\frac{\beta}{a}\right)^{(i+j)b} {_{{}_eH}}\mathcal{P}_{n-l,\beta}^{(\alpha,m,r)}\left(dx+\frac{d}{c}i, d^my, d^rz;k,a,b\right)\\
&~~\times {_{{}_eH}}\mathcal{P}_{l,\beta}^{(\alpha,m,r)}\left(cX+\frac{c}{d}j, c^mY, c^rZ;k,a,b\right)  \\
\\
=&\sum\limits_{l=0}^n d^{n-l} c^l \sum\limits_{i=0}^{d-1}\sum\limits_{j=0}^{c-1}~\left(\frac{\beta}{a}\right)^{(i+j)b} {_{{}_eH}}\mathcal{P}_{n-l,\beta}^{(\alpha,m,r)}\left(cx+\frac{c}{d}i, c^my, c^rz;k,a,b\right)\\
&~~\times {_{{}_eH}}\mathcal{P}_{l,\beta}^{(\alpha,m,r)}\left(dX+\frac{d}{c}j, d^mY, d^rZ;k,a,b\right).
\end{array}
\end{equation}
\end{thm}

\noindent
\begin{proof}Let
\begin{equation}\label{3.11}
H(t):= \frac{2^{(1-k)2\alpha}~ t^{2k\alpha}~ e^{cdxt+ y(cdt)^m} (\beta^{bc} e^{cdt}- a^{bc}) (\beta^{bd} e^{cdt} - a^{bd})~ e^{cdXt + Y(cdt)^m}}{(\beta^b e^{ct} - a^b)^{\alpha +1} ~(\beta^b e^{dt}- a^b)^{\alpha+1}~ (1-z(cdt)^r)~ (1-Z(cdt)^r)},
\end{equation}
which on rearranging the powers becomes
\begin{equation}\label{3.12}
\begin{split}
H(t) = &\frac{1}{c^{k\alpha} d^{k \alpha}}\left( \frac{2^{1-k}(ct)^k}{\beta^b e^{ct} - a^b}\right)^\alpha \frac{e^{dx(ct)+ d^my(ct)^m}}{(1-d^rz(ct)^r)} \left(\frac{\beta^{bc} e^{cdt} - a^{bc}}{\beta^b e^{dt} - a^b} \right)\\
& \times\left( \frac{2^{1-k}(dt)^k}{\beta^b e^{dt} - a^b}\right)^\alpha \frac{e^{cX(dt) + c^mY(dt)^m}}{(1-c^rZ(dt)^r)}  \left(\frac{\beta^{bd} e^{cdt} - a^{bd}}{\beta^b e^{ct} - a^b} \right).
\end{split}
\end{equation}

Now, using the series expansions for  $\left(\frac{\beta^{bc} e^{cdt} - a^{bc}}{\beta^b e^{dt} - a^b} \right)$ and $\left(\frac{\beta^{bd} e^{cdt} - a^{bd}}{\beta^b e^{ct} - a^b} \right)$  in r.h.s. of equation \eqref{3.12}, we find
\begin{equation}\label{3.13}
\begin{split}
H(t) = &\frac{1}{c^{k\alpha} d^{k \alpha}}\left( \frac{2^{1-k}(ct)^k}{\beta^b e^{ct} - a^b}\right)^\alpha \frac{e^{dx(ct)+ d^my(ct)^m}}{(1-d^rz(ct)^r)} \left(a^{(c-1)b} \sum\limits_{i=0}^{c-1} \left(\frac{\beta}{a}\right)^{ib} e^{dti} \right)\\
& \times\left( \frac{2^{1-k}(dt)^k}{\beta^b e^{dt} - a^b}\right)^\alpha \frac{e^{cX(dt) + c^mY(dt)^m}}{(1-c^rZ(dt)^r)}  \left(a^{(d-1)b} \sum\limits_{j=0}^{d-1} \left(\frac{\beta}{a}\right)^{jb} e^{ctj} \right),
\end{split}
\end{equation}
which on using equations (1.17) and then applying the Cauchy product rule gives
\begin{equation}\label{3.15}
\begin{array}{lll}
H(t) = &\frac{1}{c^{k\alpha}~ d^{k \alpha}} a^{db + bc-2b} \sum\limits_{i=0}^{c-1} \sum\limits_{j=0}^{d-1} \left(\frac{\beta}{a}\right)^{(i+j)b} \sum\limits_{n=0}^\infty \sum\limits_{l=0}^n c^{n-l} d^l {_{{}_eH}}\mathcal{P}_{n-l,\beta}^{(\alpha,m,r)}\left(dx+\frac{d}{c}i, d^my, d^rz;k,a,b\right)\\
& \times {_{{}_eH}}\mathcal{P}_{l,\beta}^{(\alpha,m,r)}\left(cX+\frac{c}{d}j, c^mY, c^rZ;k,a,b\right) t^n.
\end{array}
\end{equation}

Following the similar way, we obtain the second expansion of $H(t)$ as:
\begin{equation}\label{3.16}
\begin{array}{lll}
H(t) = &\frac{1}{d^{k\alpha}~ c^{k \alpha}} a^{cb + db-2b} \sum\limits_{i=0}^{d-1} \sum\limits_{j=0}^{c-1} \left(\frac{\beta}{a}\right)^{(i+j)b} \sum\limits_{n=0}^\infty \sum\limits_{l=0}^n d^{n-l} c^l {_{{}_eH}}\mathcal{P}_{n-l,\beta}^{(\alpha,m,r)}\left(cx+\frac{c}{d}i, c^my, c^rz;k,a,b\right)\\
& \times {_{{}_eH}}\mathcal{P}_{l,\beta}^{(\alpha,m,r)}\left(dX+\frac{d}{c}j, d^mY, d^rZ;k,a,b\right) t^n.
\end{array}
\end{equation}

Equating the coefficients of like powers of $t$ in r.h.s. of expansions \eqref{3.15} and \eqref{3.16}, we are led to
assertion \eqref{3.10}.
\end{proof}

\begin{thm} For each pair of positive integers $c,~d$ and for all integers $n \geq0, ~\alpha \geq1, ~a,b,k \in \mathbb{C}$, the following symmetry relation for the unified Apostol type-truncated exponential-Gould-Hopper polynomials ${_{{}_eH}}\mathcal{P}_{n,\beta}^{(\alpha,m,r)}(x,y,z;k,a,b)$   holds true:
\begin{equation}\label{3.17}
\begin{array}{lll}
&\hspace{-1cm}\sum\limits_{i=0}^{d-1} a^{db} \left(\frac{\beta}{a}\right)^{ib} \sum\limits_{l=0}^n {_{{}_eH}}\mathcal{P}_{n-l,\beta}^{(\alpha,m,r)}(cx, c^my, c^rz, k,a,b)d^{n-l} (ic)^l\\
=&\sum\limits_{i=0}^{c-1} a^{cb} \left(\frac{\beta}{a}\right)^{ib} \sum\limits_{l=0}^n {_{{}_eH}}\mathcal{P}_{n-l,\beta}^{(\alpha,m,r)}(dx, d^my, d^rz, k,a,b)c^{n-l} (id)^l.
\end{array}
\end{equation}
\end{thm}

\noindent
\begin{proof}Let
\begin{equation}\label{3.18}
\begin{array}{lll}
N(t):&= \frac{(2^{1-k} t^k) e^{cdxt+ y(cdt)^m} (\beta^{bd} e^{cdt} - a^{bd})}{(\beta^{b} e^{ct} - a^{b}) (\beta^{b} e^{dt} - a^{b})(1-z(cdt)^r)}\\
\\
&=\frac{1}{d^k} \left(\frac{2^{1-k} (dt)^k}{\beta^{b} e^{dt} - a^{b}}\right) \frac{e^{cx(dt)+ c^my(dt)^m}}{(1-c^rz(dt)^r)} \left(\frac{\beta^{bd} e^{cdt} - a^{bd}}{\beta^{b} e^{ct} - a^{b}}\right).
\end{array}
\end{equation}

Following the same lines of proof as in Theorem 5.2, we can prove Theorem 5.3. We omit the proof.
\end{proof}

\begin{thm} For each pair of positive integers $c,~d$ and for all integers $n \geq0, ~\alpha \geq1, ~a,b,k \in \mathbb{C}$, the following symmetry relation  between the unified Apostol type-truncated exponential-Gould-Hopper polynomials ${_{{}_eH}}\mathcal{P}_{n,\beta}^{(\alpha,m,r)}(x,y,z;k,a,b)$  involving multiple power sums holds true:
\begin{equation}\label{3.19}
\begin{array}{llll}
& \sum\limits_{l=0}^n {_{{}_eH}}\mathcal{P}_{n-l,\beta}^{(m,r)}(dx, d^my, d^rz, k,a,b)~ a^{b(d-1)\alpha} \left(\frac{\beta}{a}\right)^{- b\alpha} \sum\limits_{p=0}^l  \sum\limits_{q=0}^p {p \choose q} (-\alpha)^{p-q}~ \mathcal{S}^{(\alpha)}_k\left(d, \left(\frac{\beta}{a}\right)^b \right)\\
&~~~~~~~~~~\times {_{{}_eH}}\mathcal{P}_{l-m,\beta}^{(\alpha+1,m,r)}(cX, c^mY, c^rZ, k,a,b)~ d^{l-p-k(\alpha+1)}~ c^{n+p-(k+l)}\\
& =\sum\limits_{l=0}^n {_{{}_eH}}\mathcal{P}_{n-l,\beta}^{(m,r)}(cx, c^my, c^rz, k,a,b)~ a^{b(c-1)\alpha} \left(\frac{\beta}{a}\right)^{- b\alpha} \sum\limits_{p=0}^l  \sum\limits_{q=0}^p {p \choose q} (-\alpha)^{p-q}~ \mathcal{S}^{(\alpha)}_k\left(c, \left(\frac{\beta}{a}\right)^b \right)\\
&~~~~~~~~~~\times {_{{}_eH}}\mathcal{P}_{l-m,\beta}^{(\alpha+1,m,r)}(dX, d^mY, d^rZ, k,a,b)~ c^{l-p-k(\alpha +1)}~ d^{n+p-(k+l)}.
\end{array}
\end{equation}
\end{thm}

\noindent
\begin{proof}
Let
\begin{equation}\label{3.20}
F(t):= \frac{(2^{1-k} t^k)^{\alpha +2}~ e^{dx(ct)+ d^my (ct)^m} (\beta^{bd} e^{dct} - a^{bd})^\alpha~ e^{cX(dt)+ c^mY(dt)^m} }{(\beta^b e^{dt} -a^b)^{\alpha +1}~ (\beta^b e^{ct} - a^b)^{\alpha+1}~ (1-z(cdt)^r)(1-Z(cdt)^r)},
\end{equation}
which on rearranging the powers becomes
\begin{equation}\label{3.21}
\begin{array}{llll}
F(t):=& \frac{1}{c^k d^{k(\alpha +1)}}\left(\frac{2^{1-k} (ct)^k}{\beta^b e^{ct} - a^b}\right) \frac{e^{dx(ct)+ d^my (ct)^m}}{(1-d^rz(ct)^r)} \left(\frac{\beta^{bd} e^{dct} - a^{bd}}{\beta^b e^{ct} - a^b}\right)^\alpha \left(\frac{2^{1-k} (dt)^k}{\beta^b e^{dt} - a^b}\right)^{\alpha +1} \frac{e^{cX(dt)+ c^mY (dt)^m}}{(1-c^rZ(dt)^r)}.
\end{array}
\end{equation}

Using equations (1.17) and \eqref{3.2} in the above equation and then applying the Cauchy product rule, it follows that
\begin{equation}\label{3.23}
\begin{array}{llll}
F(t):=& \sum\limits_{n=0}^\infty \sum\limits_{l=0}^n {_{{}_eH}}\mathcal{P}_{n-l,\beta}^{(m,r)}(dx, d^my, d^rz, k,a,b) \sum\limits_{p=0}^l \sum\limits_{q=0}^p \binom{p}{q} (-\alpha)^{p-q} a^{b(d-1)\alpha} \left(\frac{\beta}{a}\right)^{- b\alpha}\\
&~~ \times \mathcal{S}_k^{(\alpha)}\left(d, \left(\frac{\beta}{a}\right)^b\right)  {_{{}_eH}}\mathcal{P}_{l-p,\beta}^{(\alpha+1,m,r)}(cX, c^mY, c^rZ, k,a,b) d^{l-p-k(\alpha+1)} \frac{c^{n+p-k-l}}{p!}.
\end{array}
\end{equation}

In a similar manner, we have
\begin{equation}\label{3.24}
\begin{array}{llll}
F(t):=& \sum\limits_{n=0}^\infty \sum\limits_{l=0}^n {_{{}_eH}}\mathcal{P}_{n-l,\beta}^{(m,r)}(cx, c^my, c^rz, k,a,b) \sum\limits_{p=0}^l \sum\limits_{q=0}^p \binom{p}{q} (-\alpha)^{p-q} a^{b(c-1)\alpha} \left(\frac{\beta}{a}\right)^{- b\alpha}\\
&~~ \times \mathcal{S}_k^{(\alpha)}\left(c, \left(\frac{\beta}{a}\right)^b\right)  {_{{}_eH}}\mathcal{P}_{l-p,\beta}^{(\alpha+1,m,r)}(dX, d^mY, d^rZ, k,a,b) c^{l-p-k(\alpha+1)} \frac{d^{n+p-k-l}}{p!}.
\end{array}
\end{equation}

Equating the coefficients of same powers of $t$ in r.h.s. of expansions \eqref{3.23} and \eqref{3.24}, we are led to assertion \eqref{3.19}.
\end{proof}

\begin{thm} For each pair of positive integers $c,~d$ and for all integers $n \geq0, ~\alpha \geq1, ~a,b,k \in \mathbb{C}$, the following symmetry relation  between the unified Apostol type-truncated exponential-Gould-Hopper polynomials ${_{{}_eH}}\mathcal{P}_{n,\beta}^{(\alpha,m,r)}(x,y,z;k,a,b)$  involving multiple power sums holds true:
\begin{equation}\label{3.25}
\begin{array}{llll}
&  \sum\limits_{p=0}^n {_{{}_eH}}\mathcal{P}_{n-p,\beta}^{(\alpha,m,r)}(dx, d^my, d^rz, k,a,b) c^{n-p-k\alpha}  \left(\frac{\beta}{a}\right)^{- b\alpha} a^{cb\alpha} \sum\limits_{q=0}^p \binom{p}{q} (-\alpha)^{p-q}~ \mathcal{S}_k^{(\alpha)}\left(c, \left(\frac{\beta}{a}\right)^b\right){d^p}\\
& =  \sum\limits_{p=0}^n {_{{}_eH}}\mathcal{P}_{n-p,\beta}^{(\alpha,m,r)}(cx, c^my, c^rz, k,a,b) d^{n-p-k\alpha}  \left(\frac{\beta}{a}\right)^{- b\alpha} a^{db\alpha} \sum\limits_{q=0}^p \binom{p}{q} (-\alpha)^{p-q}~ \mathcal{S}_k^{(\alpha)}\left(d, \left(\frac{\beta}{a}\right)^b\right){c^p}.
\end{array}
\end{equation}
\end{thm}

\noindent
\begin{proof}
Let
\begin{equation}\label{3.26}
\begin{array}{llll}
M(t)&:=\frac{(2^{1-k} t^k)^\alpha~e^{dx(ct) +d^my(ct)^m} (\beta^{bc} e^{cdt} - a^{bc})^\alpha}{(\beta^b e^{dt}- a^b)^\alpha~ (\beta^b e^{ct}-a^b)^\alpha~ (1-z(dct)^r)}\\
\\
& = \frac{1}{c^{k\alpha}}\left( \frac{2^{1-k} (ct)^k}{\beta^b e^{ct}-a^b}\right)^\alpha \frac{e^{dx(ct) +d^my(ct)^m}}{(1-d^rz(ct)^r)} a^{(c-1)b\alpha} \left( \frac{\left(\frac{\beta}{a}\right)^{bc} e^{cdt} -1}{\left(\frac{\beta}{a}\right)^{b} e^{dt} -1}\right)^\alpha.
\end{array}
\end{equation}

Following the same lines of proof as in Theorem 3.4, we can prove Theorem 3.5. We omit the proof.
\end{proof}

\begin{thm} For each pair of positive integers $c,~d$ and for all integers $n \geq0, ~\alpha \geq1, ~a,b,k \in \mathbb{C}$, the following symmetry relation  between the unified Apostol type-truncated exponential-Gould-Hopper polynomials ${_{{}_eH}}\mathcal{P}_{n,\beta}^{(\alpha,m,r)}(x,y,z;k,a,b)$  involving Hurwitz-Lerch Zeta function $\Phi_\mu(z,s,a)$ holds true:
\begin{equation}\label{a}
\begin{array}{llll}
& \sum\limits_{n=0}^h \sum\limits_{s=0}^{h-n} \binom{h-n}{s}~ \Phi_\alpha\left( \left(\frac{\beta}{a} \right)^b, s-h+n, cx \right) H_s^{(m)}(0, c^my) ~  \left(\frac{\beta}{a}\right)^{- b}a^{(c-1)b} \\
& \sum\limits_{l=0}^n \sum\limits_{q=0}^{n-l} \binom{n-l}{q} (-1)^{n-l-q}~  \mathcal{S}_k\left(c, \left(\frac{\beta}{a}\right)^b \right) {_{{}_eH}}\mathcal{P}_{l,\beta}^{(\alpha,m,r)}(dX, d^mY, d^rZ, k,a,b)~ c^{l-k\alpha} ~d^{h-l}\\
&= \sum\limits_{n=0}^h \sum\limits_{s=0}^{h-n} \binom{h-n}{s}~ \Phi_\alpha\left( \left(\frac{\beta}{a} \right)^b, s-h+n, dx \right) H_s^{(m)}(0, d^my)  ~  \left(\frac{\beta}{a}\right)^{- b} a^{(d-1)b} \\
&\sum\limits_{l=0}^n \sum\limits_{q=0}^{n-l} \binom{n-l}{q} (-1)^{n-l-q}~  \mathcal{S}_k\left(d, \left(\frac{\beta}{a}\right)^b \right) {_{{}_eH}}\mathcal{P}_{l,\beta}^{(\alpha,m,r)}(cX, c^mY, c^rZ, k,a,b)~ d^{l-k\alpha}~ c^{h-l}.
\end{array}
\end{equation}
\end{thm}

\noindent
\begin{proof}
Let
\begin{equation}\label{3.27}
P(t):= \frac{(2^{1-k} t^k)^{2\alpha}~ e^{cx(dt)+c^my(dt)^m} (\beta^{bc} e^{cdt}-a^{bc})~ e^{dX(ct)+ d^mY(ct)^m}}{(\beta^b e^{dt} - a^b)^{\alpha+1}~(\beta^b e^{ct} -a^b)^\alpha~(1-Z(cdt)^r)},
\end{equation}
which on rearranging the powers becomes
\begin{equation}\label{3.28}
\begin{array}{llll}
P(t):=& \frac{1}{c^{k\alpha}}(2^{1-k} t^k)^\alpha~ (\beta^b e^{dt} -a^b)^{-\alpha}~ e^{cx(dt)+c^my(dt)^m} \left(\frac{\beta^{bc} e^{cdt}-a^{bc}}{\beta^{b} e^{dt}-a^{b}}\right)\\
& \times \left( \frac{2^{1-k} (ct)^k}{\beta^b e^{ct} -a^b}\right)^\alpha \frac{e^{dX(ct)+ d^mY(ct)^m}}{(1-d^r Z(ct)^r)}.
\end{array}
\end{equation}
Now, using equations (1.17), \eqref{3.2} (for $\alpha=1$) in equation (5.25) and formulae
\begin{equation}\label{3.29}
(1+w)^{-\alpha}=\sum\limits_{m=0}^\infty {m+\alpha-1 \choose m} (-w)^m;~~~~|w|<1,
\end{equation}
in the resultant equation, we find
\begin{equation}\label{3.30}
\begin{array}{llll}
P(t):=& \frac{1}{c^{k \alpha}}~(2^{1-k} t^k)^\alpha~(-a^b)^{-\alpha} \sum\limits_{p=0}^\infty \binom{p+\alpha -1}{p} \left(\frac{\beta}{a}\right)^{bp}~ e^{pdt}~e^{cx(dt)+c^my(dt)^m} a^{(c-1)b} \left(\frac{\beta}{a}\right)^{- b}  \\
& \times \sum\limits_{n=0}^\infty \sum\limits_{q=0}^n \binom{n}{q} (-1)^{n-q}~ \mathcal{S}_k\left(c, \left(\frac{\beta}{a}\right)^b \right) \frac{(dt)^n}{n !} \sum\limits_{l=0}^\infty {_{{}_eH}}\mathcal{P}_{l,\beta}^{(\alpha,m,r)}(dX, d^mY, d^rZ, k,a,b) (ct)^l.
\end{array}
\end{equation}

Making use of equation \eqref{1.11} in the above equation and then applying the Cauchy product rule, we find
\begin{equation}\label{3.32}
\begin{array}{llll}
&P(t):= \left(\frac{2^{1-k} t^k}{c^k}\right)^\alpha(-a^b)^{-\alpha} \sum\limits_{h=0}^\infty \left(\sum\limits_{s=0}^h  \binom{h}{s} \Phi_\alpha\left( \left(\frac{\beta}{a} \right)^b, s-h, cx \right) H_s^{(m)}(0, c^my) d^h \right) \frac{t^h}{h!}~ a^{(c-1)b} \\
&\hspace{-.5cm}\times \left(\frac{\beta}{a}\right)^{- b} \sum\limits_{n=0}^\infty \left(\sum\limits_{l=0}^n \sum\limits_{q=0}^{n-l} \binom{n-l}{q} (-1)^{n-l-q}  \mathcal{S}_k\left(c, \left(\frac{\beta}{a}\right)^b \right) {_{{}_eH}}\mathcal{P}_{l,\beta}^{(\alpha,m,r)}(dX, d^mY, d^rZ, k,a,b) c^l \frac{d^{n-l}}{(n-l)!} \right) t^n,
\end{array}
\end{equation}
which on again applying the Cauchy product rule gives
\begin{equation}\label{3.33}
\begin{array}{llll}
&\hspace{-.4cm}P(t):= (2^{1-k} t^k)^\alpha~(-a^b)^{-\alpha} \sum\limits_{h=0}^\infty \Bigg(\sum\limits_{n=0}^h \sum\limits_{s=0}^{h-n} \binom{h-n}{s} \Phi_\alpha\left( \left(\frac{\beta}{a} \right)^b, s-h+n, cx \right) H_s^{(m)}(0, c^my) \left(\frac{\beta}{a}\right)^{- b}\\
&\hspace{-.5cm} \times a^{(c-1)b} \sum\limits_{l=0}^n \sum\limits_{q=0}^{n-l} \binom{n-l}{q} (-1)^{n-l-q}  \mathcal{S}_k\left(c, \left(\frac{\beta}{a}\right)^b \right) {_{{}_eH}}\mathcal{P}_{l,\beta}^{(\alpha,m,r)}(dX, d^mY, d^rZ, k,a,b) c^{l-k\alpha} \frac{d^{h-l}}{(n-l)!} \Bigg) \frac{t^h}{(h-n)!}.
\end{array}
\end{equation}

In a similar manner, we have
\begin{equation}\label{3.34}
\begin{array}{llll}
&\hspace{-.4cm}P(t):= (2^{1-k} t^k)^\alpha~(-a^b)^{-\alpha} \sum\limits_{h=0}^\infty \Bigg(\sum\limits_{n=0}^h \sum\limits_{s=0}^{h-n} \binom{h-n}{s} \Phi_\alpha\left( \left(\frac{\beta}{a} \right)^b, s-h+n, dx \right) H_s^{(m)}(0, d^my)  ~ \left(\frac{\beta}{a}\right)^{- b}\\
&\hspace{-.5cm} \times a^{(d-1)b} \sum\limits_{l=0}^n \sum\limits_{q=0}^{n-l} \binom{n-l}{q} (-1)^{n-l-q}  \mathcal{S}_k\left(d, \left(\frac{\beta}{a}\right)^b \right) {_{{}_eH}}\mathcal{P}_{l,\beta}^{(\alpha,m,r)}(cX, c^mY, c^rZ, k,a,b) d^{l-k\alpha} \frac{c^{h-l}}{(n-l)!} \Bigg) \frac{t^h}{(h-n)!}.
\end{array}
\end{equation}

Equating the coefficients of same powers of $t$ in r.h.s. of expansions \eqref{3.33} and \eqref{3.34}, we are led to assertion \eqref{a}.\\
\end{proof}

In view of Remarks 1.2-1.4, the symmetry identities for the TEGHABP ${_{{}_eH}}\mathfrak{B}_{n}^{(\alpha,m,r)}(x,y,z;\lambda)$, TEGHAEP ${_{{}_eH}}\mathfrak{E}_{n}^{(\alpha,m,r)}(x,y,z;\lambda)$  and TEGHAGP ${_{{}_eH}}\mathcal{G}_{n}^{(\alpha,m,r)}(x,y,z;\lambda)$ can be obtained by suitably substituting the values of the parameters $k,a,b$ and $\beta$.

\section{Conclusion}
The formulas and identities provided in this paper confirm the usefulness of the series rearrangement techniques. We have derived several unified formulas and identities for the unified Apostol type-truncated exponential-Gould-Hopper polynomials ${_{{}_eH}}\mathcal{P}_{n,\beta}^{(\alpha,m,r)}(x,y,z;k,a,b)$  by using different analytical means on their generating functions. This process could be used to establish further quite a wide variety of formulas for several other special polynomials and can be extended to derive new relations for conventional and generalized polynomials.

\end{document}